\documentstyle[12pt,leqno,amsfonts,amscd]{amsart}
\newtheorem{Theorem}{Theorem}[section]

\newtheorem{Proposition}[Theorem]{Proposition}

\newtheorem{Lemma}[Theorem]{Lemma}
\newtheorem{Corollary}[Theorem]{Corollary}
\theoremstyle{remark}

\def\diam{\operatorname{diam}}
\begin{document}
\title{Overinterpolation}
\author{ Dan Coman and Evgeny A. Poletsky}
\thanks{Both authors are supported by NSF Grants.}
\subjclass[2000]{Primary: 30E05. Secondary: 30E10, 41A25}
\address{ Department of Mathematics,  215 Carnegie Hall,
Syracuse University,  Syracuse, NY 13244-1150, USA. E-mail:
dcoman@@syr.edu, eapolets@@syr.edu}
\begin{abstract} In this paper we study the consequences of
overinterpolation, i.e., the situation when a function can be
interpolated by polynomial, or rational, or algebraic functions in
more points that normally expected. We show that in many cases
such a function has specific forms.
\end{abstract}
\maketitle

\section{Introduction}
\par Let ${\mathcal P}_n$ be the space of all polynomials on the
complex plane ${\Bbb C}$ whose degree is at most $n$. Let
${\mathcal R}_{nm}$ be the space of rational functions
$R_{nm}=P_n/Q_m$ where $P_n\in{\mathcal P}_n$ and $Q_m\in{\mathcal
P}_m$.
\par If $f$ is a function on a compact set $K\subset{\Bbb C}$,
then we denote by $N_K(n)$ and $N_K(n,m)$ the maximal number of
zeros on $K$ of the functions $f-p$, where $p\in{\mathcal P}_n$,
respectively $p\in{\mathcal R}_{nm}$. Since functions in
${\mathcal P}_n$ or ${\mathcal R}_{nm}$ have $n+1$ or,
respectively, $n+m+2$ coefficients, $N_K(n)\ge n+1$ and
$N_K(n,m)\ge n+1$.
\par In this paper we consider the situations when for a fixed
function $f$ we have either polynomial or rational {\it
overinterpolation}. This means that
$\lim_{n\to\infty}N_K(n)/n=\infty$ or
$\lim_{n\to\infty}N_K(n,m)/n=\infty$.
\par One can expect that in the case of overinterpolation, the
function $f$ must be either polynomial  or rational. We prove two
theorems of this kind. Before we state them, let us introduce some
notation.
\par Let $\Delta_r\subset{\Bbb C}$ be the open disk of radius $r$
centered at the origin, and $\Delta\subset{\Bbb C}$ be the open
unit disk. We denote by $O(\Delta_r)$ and $O(\overline\Delta_r)$
the set of holomorphic functions on $\Delta_r$, respectively on
neighborhoods of $\overline\Delta_r$.
\par The first theorem proved in Section \ref{S:obp} states that
for an analytic function $f$ overinterpolation by polynomials
implies that $f$ is a polynomial.
\begin{Theorem}\label{T:Lagr}Let $f\in O(\Delta)$ and
$K=\overline\Delta_r$, where $r<1$. If
$\lim_{n\to\infty}N_K(n)/n=\infty$ then $f$ is a polynomial.
\end{Theorem}
\par For a function $f$ as above, the second theorem states that
overinterpolation in ${\mathcal R}_{n1}$ implies that either $f$
is entire or it belongs to ${\mathcal R}_{n1}$.
\begin{Theorem}\label{T:rational} Let $f\in O(\Delta)$ and
$K=\overline\Delta_r$, where $r<1$. If
$\lim_{n\to\infty}N_K(n,1)/n=\infty$ then either $f$ is entire or
$f=P/Q$, where $P,\,Q$ are polynomials, $\deg Q=1$ and $Q$ does
not divide $P$.
\end{Theorem}
\par This theorem is proved in Section \ref{S:obrf}, where we also
consider the case of Pad\'e interpolation, i.e. when $K=\{0\}$.
For any germ of an analytic function $f$ at 0 and any fixed
$m\in{\Bbb N}$, we show that $N_K(n,m)\leq n+m+1$ for infinitely
many $n$, unless $f$ is the germ of a rational function.
\par The expected rate of rational approximation of continuous or
analytic functions is at most geometric, but in some cases
functions can be approximated faster. This phenomenon is called
{\it overconvergence}. In \cite{Go} and \cite{Ch} Gonchar and
Chirka have shown that in this case the functions have specific
forms. In Section \ref{S:oao} we prove that overinterpolation
implies overconvergence on some circle and, therefore,
overinterpolated functions have the same specific forms as in the
results of Gonchar and Chirka.
\par For the entire function $f(z)=\sum 2^{-n!}z^n$, its Taylor
series is overconvergent but by Theorem \ref{T:Lagr} $f$ cannot be
overinterpolated by polynomials. Hence overconvergence does not
imply overinterpolation.
\par The assumption in all our results that $f\in O(\Delta)$ seems
to be a technical necessity. In the last section we consider the
interpolation of a general set $S$ in ${\Bbb C}^2$ by algebraic
functions, i. e., we are looking for the maximal number $N_S^a(n)$
of zeros on $S$ of a polynomial of degree $n$ which does not
vanish on $S$. The desirable estimate is $N_S^a(n)\le An^\alpha$,
where $A$ and $\alpha$ are some constants. We show that either $S$
is finite, or $\alpha=1$ and $S$ is contained in an irreducible
algebraic curve, or $\alpha\ge2$.
\par It should be noted that in \cite{CP3} we proved for a
large class of meromorphic functions $f$ on ${\Bbb C}$ with
finitely many poles, including the Riemann $\zeta$-function, that
if $S$ is the graph of $f$ over $\Delta_r$, then $N^a_S(n)\le
An^2\log r$.

\section{Overinterpolation by polynomials}\label{S:obp}
If $f\in O(\overline\Delta_R)$ we set
$$M(r,f)=\max\{|f(z)|:\,|z|=r\},\;r\leq R.$$

\par We will need the following lemmas:
\begin{Lemma}\label{L:Lagr} Let $f\in O(\overline\Delta_R)$ and
$L_nf$ denote the Lagrange interpolating polynomial of $f$ at the
(not necessarily distinct) points
$z_0,\dots,z_n\in\overline\Delta_r$, where $r<R$. If $0<s<R$ then
$$M(s,f-L_nf)\leq
M(R,f)\frac{R}{R-s}\left(\frac{s+r}{R-r}\right)^{n+1}.$$\end{Lemma}
\begin{pf} Let $\omega(z)=(z-z_0)\dots(z-z_n)$. By \cite[p. 59, (1.4)]{G}
we have $$f(z)-L_nf(z)=\frac{1}{2\pi i}\int_{|t|=R}
\frac{\omega(z)f(t)}{\omega(t)(t-z)}\,dt.$$ The lemma follows
since $|\omega(t)|\geq(R-r)^{n+1}$ for $|t|=R$, and since
$M(s,\omega)\leq(s+r)^{n+1}$.\end{pf}
\par For $R>0$ let $$R_+=\max\{R,1\}.$$ We have the following estimate of
Taylor coefficients.
\begin{Lemma}\label{L:eftc} Let $f\in O(\Delta)$,
$f(z)=\sum_{k\ge0}f_kz^k$. Suppose that $0<r<1$ and the function
$f-P_n$ has $N$ zeros in $\overline\Delta_r$, where $P_n$ is a
polynomial of degree at most $n$. There exist positive constants
$A\ge 1$, $a<1$ and $\delta$, depending only on $r$, with the
following property: If $N\ge A(n+1)$, then
$$|f_k|\le \frac{M(R,f)}{R_+^{n+1}}\;a^N,$$ for $n<k\le\delta N$ and
every $R\ge(r+2)/3$ such that $f\in O(\overline\Delta_R)$.
\end{Lemma}
\begin{pf} Let $s=(2r+1)/3$ and fix $R$ as in the
statement. Let $z_0,\dots,z_{N-1}$ be zeros of $f-P_n$ in
$\overline\Delta_r$. Since $N\ge n+1$ the polynomial $P_n=L_nf$ is
the Lagrange interpolating polynomial of $f$ at $z_0,\dots,z_n$.
Since $f-P_n$ has $N$ zeros in $\overline\Delta_r$, we have by
\cite[Theorem 2.2]{CP2} (see the formula on p. 578)
$$M(r,f-P_n)\leq
M(s,f-P_n)\left(\frac{2rs}{r^2+s^2}\right)^{N}.$$ Hence by Lemma
\ref{L:Lagr}
\begin{eqnarray*}M(r,f-P_n)&\leq&
M(R,f)\frac{R}{R-s}\left(\frac{s+r}{R-r}\right)^{n+1}
\left(\frac{2rs}{r^2+s^2}\right)^{N}\\
&=&\frac{M(R,f)}{R^{n+1}}\frac{1}{1-s/R}
\left(\frac{s+r}{1-r/R}\right)^{n+1}
\left(\frac{2rs}{r^2+s^2}\right)^{N}. \end{eqnarray*} Notice that
$$\frac{1}{1-s/R}<\frac3{1-r}\;,\;
\frac{s+r}{1-r/R}<\frac3{1-r}\;,$$ and
$$a_1:=\frac{2rs}{r^2+s^2}<1.$$
Since $R>2/3$ we obtain
\begin{eqnarray}\label{e:mfe0}M(r,f-P_n)&\leq&\frac{M(R,f)}{R^{n+1}}
\left(\frac{3}{1-r}\right)^{n+2}a^N_1\\
&<&\frac{M(R,f)}{R_+^{n+1}}\left(\frac{3}{2}\right)^{n+1}
\left(\frac3{1-r}\right)^{n+2}a^N_1\notag\\
&<&\frac{M(R,f)}{R_+^{n+1}}\left(\frac{5}{1-r}\right)^{n+2}a^N_1.
\notag\end{eqnarray} Let $$A=\max\left\{-4\,\frac{\log
5-\log(1-r)}{\log a_1}\,,1\right\}\;.$$ As $N\ge A(n+1)$ we obtain
$$M(r,f-P_n)\leq\frac{M(R,f)}{R_+^{n+1}}\;a_2^N,$$
where $a_2=a_1^{1/2}$. Since $k>\deg P_n$ it follows by Cauchy's
inequalities that
$$|f_k|\leq\frac{M(r,f-P_n)}{r^k}\leq
\frac{M(R,f)}{R_+^{n+1}}\;a_2^Nr^{-k}.$$ We define $a=a^{1/2}_2$
and $\delta$ by $r^\delta=a$. If $k\le\delta N$ then
$$|f_k|\leq\frac{M(R,f)}{R_+^{n+1}}\;a^N.$$
\end{pf}

\vspace{2mm}\noindent{\em Proof of Theorem \ref{T:Lagr}.} Let
$f(z)=\sum_{n\geq0}f_nz^n$. We can find an increasing sequence of
integers $N(n)\leq N_K(n)$ such that $N(n)/n\rightarrow\infty$ and
the function $f-P_n$ has at least $N(n)$ zeros in
$\overline\Delta_r$, where $P_n\in{\mathcal P}_n$.
\par Fix $R\geq(r+2)/3$ so that $f\in O(\overline\Delta_R)$. Let
$a<1\le A$ be the constants from Lemma \ref{L:eftc} and $n_0$ be
so that $N(n)\geq A(n+1)$ if $n\geq n_0$. Lemma \ref{L:eftc}
implies that for $n\geq n_0$
\begin{equation}\label{E:Lagr}
|f_{n+1}|\le \frac{M(R,f)}{R_+^{n+1}}\;a^{N(n)}.\end{equation}
Therefore $|f_n|^{1/n}\rightarrow0$, so $f$ is entire, hence
(\ref{E:Lagr}) holds for any $R\geq1$. By Cauchy's inequalities we
have $|f_n|\leq M(R,f)/R^n$ for $n\leq n_0$. Using these estimates
of the coefficients, we obtain the following bound for $M(2R,f)$,
$R\geq1$:
\begin{equation}\label{E:doubling}M(2R,f)\leq\sum_{n\geq0}
|f_n|(2R)^n\leq CM(R,f),\end{equation} where
$$C=\sum_{n=0}^{n_0}2^n+\sum_{n=n_0}^\infty 2^{n+1}a^{N(n)}$$
is independent on $R$. Note that
$$2^na^{N(n)}=\left(2a^{N(n)/n}\right)^n\leq2^{-n},$$ provided
that $n$ is sufficiently large, thus $C$ is finite. \par Applying
the doubling inequality (\ref{E:doubling}) successively we obtain
$$M(2^j,f)\leq C^jM(1,f),$$ for any $j>0$. Hence $$|f_n|\leq
\frac{C^jM(1,f)}{2^{nj}}\rightarrow0\;{\rm
as}\;j\rightarrow\infty,$$ provided that $2^n>C$. We conclude that
$f$ is a polynomial of degree at most $\log C/\log 2$. $\Box$
\par Theorem \ref{T:Lagr} has the following immediate corollary:
\begin{Corollary}\label{C:Lagr} Let $\{n_k\}_{k\geq0}$ be an
increasing sequence of natural numbers such that $n_{k+1}/n_k\leq
C$ for some constant $C$. Let $f\in O(\Delta)$ and
$K=\overline\Delta_r$, where $r<1$. If
$\lim_{k\to\infty}N_K(n_k)/n_k=\infty$ then $f$ is a
polynomial.\end{Corollary}
\begin{pf} Let $N(n)=N_K(n)$. If $n_k\leq n<n_{k+1}$ then
$$\frac{N(n)}{n}>\frac{N(n_k)}{n_{k+1}}\geq\frac{N(n_k)}{Cn_k}\;.$$
\end{pf}

\section{Overinterpolation by rational functions}\label{S:obrf}
\par We prove here Theorem \ref{T:rational}. We can find an
increasing sequence of integers $N(n)\leq N_K(n,1)$ such that
$N(n)/n\rightarrow\infty$ and the function $Q_nf-P_n$ has at least
$N(n)$ zeros in $\overline\Delta_r$, where $P_n\in{\mathcal P}_n$,
$Q_n\in{\mathcal P}_1$ and $Q_n\neq0$. \par Let us write
$$f(z)=\sum_{k\geq0}f_kz^k,\;Q_n(z)=\alpha_nz-\beta_n.$$
Let $$\rho=\frac{1}{\limsup|f_k|^{1/k}}\geq1$$ be the radius of
convergence of the power series of $f$ at the origin. \par By
considering functions $c(Q_nf-P_n)$, where $c\in{\Bbb
C}\setminus\{0\}$, we can identify $Q_n$ with the point
$[\alpha_n:\beta_n]\in{\Bbb P}^1$. Thus
$$Q_n(z)=\alpha_nz-1,\;\alpha_n\in{\Bbb C},\; {\rm
or}\;Q_n(z)=z.$$ The latter case corresponds to $\alpha_n=\infty$
in the extended complex plane.
\par We begin with a few lemmas.
\begin{Lemma}\label{L:eftc2} There exist constants
$a<1$, $\delta<1$ and an integer $n_0$, depending only on $r$,
with the following property: If $n\geq n_0$, then one of the
inequalities
$$|\alpha_nf_{k-1}-f_k|\leq\frac{M(R,f)}{R_+^{n+1}}\;
(|\alpha_n|R_++1)a^{N(n)}\;,\;\;|f_{k-1}|\leq
\frac{M(R,f)}{R_+^n}\;a^{N(n)},$$ holds for every $k$,
$n<k\leq\delta N(n)$, and every $R\geq(r+2)/3$ such that $f\in
O(\overline\Delta_R)$.
\end{Lemma}
\begin{pf} Let $a<1\le A$, $\delta>0$, be the constants from Lemma
\ref{L:eftc}, and let $n_0=n_0(r)$ be an integer such that
$N(n)\geq A(n+1)$ for $n\geq n_0$. We fix such an $n$, and apply
Lemma \ref{L:eftc} to the function $Q_nf$ and the polynomial
$P_n$. If $Q_n(z)=\alpha_nz-1$ then
$$Q_n(z)f(z)=-f_0+\sum_{k\geq1}(\alpha_nf_{k-1}-f_k)z^k,$$ and
$M(R,Q_nf)\leq M(R,f)(|\alpha_n|R_++1)$. This yields the first
inequality of the lemma. The second one is obtained in a similar
way, in the case when $Q_n(z)=z$ (or by letting
$\alpha_n\rightarrow\infty$).
\end{pf}
\begin{Lemma}\label{L:liminf} If $f\in O(\Delta_s)$, $1\leq
s\leq\infty$, and if $\;\liminf_{n\rightarrow\infty}
|\alpha_n|>1/s$, then $f$ is a polynomial.\end{Lemma}
\begin{pf} There exist $n_1\geq n_0$ and $\epsilon>0$ such that
$|\alpha_n|>1/s+\epsilon$, for $n\geq n_1$. Let
$$c=\left(\frac{1}{s}+\epsilon\right)^{-1},\;d=1+c.$$
By Lemma \ref{L:eftc2} with $k=n+1$ we have
$$|f_n|\leq\frac{|f_{n+1}|}{|\alpha_n|}+
\frac{M(R,f)}{R_+^n}\left(1+\frac{1}{R_+|\alpha_n|}\right)a^{N(n)}
\leq c|f_{n+1}|+\frac{d M(R,f)}{R_+^n}\;a^{N(n)},$$ for every
$R\geq(r+2)/3$ such that $f\in O(\overline\Delta_R)$. Note that
this estimate obviously holds in the case $\alpha_n=\infty$.
Applying it successively we obtain
\begin{equation}\label{E:liminf}|f_n|\leq c^k|f_{n+k}|+\frac{d
M(R,f)}{R_+^n}\;\sum_{j=0}^{k-1}\frac{c^j}{R_+^j}\;a^{N(n+j)},
\end{equation} for every $k\geq1$. Fix $s_1\geq(r+2)/3$ such
that $c<s_1<s$. Since $f\in O(\Delta_s)$ we have
$$|f_{n+k}|\leq\left(\frac{1}{s}+\frac{\epsilon}{2}\right)^{n+k},$$
for $k$ sufficiently large. Since $N(n)$ is increasing, and if
$R\geq s_1$, we obtain by (\ref{E:liminf}) $$|f_n|\leq
c^k|f_{n+k}|+\frac{d
M(R,f)}{R_+^n}\;a^{N(n)}\sum_{j=0}^\infty\frac{c^j}{s_1^j}
\leq\frac{\left(\frac{1}{s}+\frac{\epsilon}{2}\right)^{n+k}}
{\left(\frac{1}{s}+\epsilon\right)^k}+\frac{ds_1
M(R,f)}{(s_1-c)R_+^n}\;a^{N(n)}.$$ Letting $k\rightarrow\infty$ we
conclude that
$$|f_n|\leq\frac{ds_1M(R,f)}{(s_1-c)R_+^n}\;a^{N(n)}$$
holds for all $n\geq n_1$ and $R\geq s_1$ such that $f\in
O(\overline\Delta_R)$. This is a similar estimate to
(\ref{E:Lagr}) from the proof of Theorem \ref{T:Lagr}. Therefore,
by the same argument as in the proof of Theorem \ref{T:Lagr}, it
follows that $f$ is entire and $M(2R,f)\leq CM(R,f)$ for every
$R\geq s_1$, where
$$C=\sum_{n=0}^{n_1-1}2^n+\frac{ds_1}{s_1-c}\; \sum_{n=n_1}^\infty
2^na^{N(n)}$$ is independent on $R$. Hence $f$ is a
polynomial.\end{pf}
\begin{Lemma}\label{L:limsup} $\limsup_{n\rightarrow\infty}
|\alpha_n|\geq1/\rho$.\end{Lemma}
\begin{pf} We assume for a contradiction that there exist $n_1\geq
n_0$ and $0<\epsilon<1/\rho$ such that
$$|\alpha_n|<c:=\rho^{-1}-\epsilon,\;n\geq n_1.$$ As $c<1$,
we obtain by Lemma \ref{L:eftc2}, applied with $k=n+1$ and
$R=(r+2)/3<1$, that
$$|f_{n+1}|\leq|\alpha_nf_n|+M(|\alpha_n|+1)a^{N(n)}\leq
c|f_n|+2Ma^{N(n)},$$ where $M=M(R,f)$ and $n\geq n_1$. Hence
$$|f_{n+k}|\leq c^k|f_n|+2M\sum_{j=0}^{k-1}c^ja^{N(n+k-1-j)}
\leq c^k|f_n|+2Ma^{N(n)}\sum_{j=0}^\infty c^j,$$ for all $k\geq1$.
\par Let $C=2M/(1-c)$. Then for $n\geq n_1$ we have
$$|f_{2n}|\leq c^n|f_n|+Ca^{N(n)}\;,\;\;|f_{2n+1}|\leq
c^{n+1}|f_n|+Ca^{N(n)}.$$ Since, for $n$ large,
$|f_n|\leq(\rho^{-1}+\epsilon)^n$, it follows that
\begin{eqnarray*}
|f_{2n}|^{1/(2n)}&\leq&c^{1/2}(\rho^{-1}+\epsilon)^{1/2}+
C^{1/(2n)}a^{N(n)/(2n)},\\|f_{2n+1}|^{1/(2n+1)}&\leq&
c^{(n+1)/(2n+1)}(\rho^{-1}+\epsilon)^{n/(2n+1)}+
C^{1/(2n+1)}a^{N(n)/(2n+1)}.\end{eqnarray*} Note that
$a^{N(n)/n}\rightarrow0$. Therefore
$$\rho^{-1}=\limsup_{j\rightarrow\infty}|f_j|^{1/j}\leq
(\rho^{-1}-\epsilon)^{1/2}(\rho^{-1}+\epsilon)^{1/2},$$ a
contradiction.\end{pf}

\vspace{2mm}\noindent{\em Proof of Theorem \ref{T:rational}.} We
can assume $\rho<\infty$, otherwise $f$ is entire. The radius of
convergence of the power series of $f(\rho z)$ at the origin is 1,
and the function $Q_n(\rho z)f(\rho z)-P_n(\rho z)$ has $N(n)$
zeros in the disk $\overline\Delta_{r/\rho}$. Therefore we may
assume that $\rho=1$.
\par Let $R=(r+2)/3$ and $M=M(R,f)$. By Lemma \ref{L:eftc2}, one
of the estimates
\begin{equation}\label{E:mainest}|\alpha_nf_{k-1}-f_k|\leq
M(|\alpha_n|+1)a^{N(n)}, \;|f_{k-1}|\le Ma^{N(n)},\end{equation}
holds for $n<k\leq\delta N(n)$, provided that $n\geq n_0$.
\par By Lemma \ref{L:limsup}, $|\alpha_n|>1/3$ or $\alpha_n=\infty$
for infinitely many $n$. We show that there exists a sequence
$m_j\to\infty$ such that $$\alpha_{m_j}\in{\Bbb C},\;
|\alpha_{m_j}|>1/3,\;|f_{m_j+1}|>2^{-m_j-2}.$$ Fix any $n$ large
with $|\alpha_n|>1/3$ or $\alpha_n=\infty$. Let $k\ge n$ be the
smallest integer such that $|f_k|>2^{-k}$. Such $k$ exists since
$\rho=1$. If $k>n$ then $\alpha_{k-1}$ is finite. Otherwise by
(\ref{E:mainest})
$$2^{-k}<|f_k|\le Ma^{N(k-1)},$$ which is impossible as $n$ is
large. By the definition of $k$, $|f_{k-1}|\leq2^{-k+1}$. We claim
that $|\alpha_{k-1}|>1/3$. If not, then using (\ref{E:mainest})
$$2^{-k}<|f_k|\le|\alpha_{k-1}f_{k-1}|+2Ma^{N(k-1)}\le
\frac{2^{-k+1}}{3}+2Ma^{N(k-1)},$$ so $2^{-k}<6Ma^{N(k-1)}$. This
is a contradiction since $n$ is large.
\par If $k=n$, then $|f_n|>2^{-n}>Ma^{N(n)}$ shows that
$\alpha_n\in{\Bbb C}$, so $|\alpha_n|>1/3$. We have by
(\ref{E:mainest})
\begin{eqnarray*}|f_{n+1}|&\ge&|\alpha_nf_n|-4M|\alpha_n|a^{N(n)}>
|\alpha_n|(2^{-n}-4Ma^{N(n)})\\
&>&\frac{2^{-n}}{3}\left(1-2^{n+2}Ma^{N(n)}\right)\ge
2^{-n-2}.\end{eqnarray*} This establishes the existence of the
desired sequence $m_j$.
\par Since $N(n)/n\rightarrow\infty$ and $\rho=1$, we can find
$n_1\geq n_0$ with the property that
\begin{equation}\label{E:ass}
2^{3n+6k+16}Ma^{N(n+k)}<1,\;|f_n|<2^n,
\end{equation}
hold for every $n\ge n_1$ and for every $k\ge0$. Then we fix $n\ge
n_1$ such that $\alpha_n\in{\Bbb C}$, $|\alpha_n|>1/3$ and
$|f_{n+1}|>2^{-n-2}$. We have using (\ref{E:mainest}) that
$$2^{-n-2}|\alpha_n|<|\alpha_{n}||f_{n+1}|\le|f_{n+2}|+
4M|\alpha_n|a^{N(n)}\leq2^{n+2}+4M|\alpha_n|a^{N(n)},$$ so by
(\ref{E:ass})
$$|\alpha_{n}|\leq\frac{2^{2n+4}}{1-2^{n+4}Ma^{N(n)}}\leq2^{2n+5}.$$
\par We will show by induction that for every $k\geq0$
\begin{equation}\label{E:ind}|f_{n+k+1}|>2^{-n-2}6^{-k}\;,\;\;
\frac16+6^{-k-1}<|\alpha_{n+k}|<2^{2n+6}-6^{-k-1}.\end{equation}
Evidently, these inequalities hold for $k=0$. Suppose that they
are true for some $k\geq0$. Then using (\ref{E:mainest})
\begin{eqnarray*}|f_{n+k+2}|&\geq&|\alpha_{n+k}f_{n+k+1}|-
7M|\alpha_{n+k}|a^{N(n+k)}\\&\geq&
2^{-n-2}6^{-k-1}+2^{-n-2}6^{-2k-1}-2^{2n+9}Ma^{N(n+k)}.\end{eqnarray*}
By (\ref{E:ass})
$$2^{-n-2}6^{-2k-1}-2^{2n+9}Ma^{N(n+k)}>0,$$ so we see that
$|f_{n+k+2}|>2^{-n-2}6^{-k-1}$.
\par Since $|f_{n+k+1}|>2^{-n-2}6^{-k}$, we have in view of
(\ref{E:mainest}) and (\ref{E:ass}) that $\alpha_{n+k+1}\in{\Bbb
C}$. Therefore by (\ref{E:mainest})
$$|\alpha_{n+k}f_{n+k+1}-f_{n+k+2}|\le
M(|\alpha_{n+k}|+1)a^{N(n+k)},$$
$$|\alpha_{n+k+1}f_{n+k+1}-f_{n+k+2}|\le
M(|\alpha_{n+k+1}|+1)a^{N(n+k+1)}.$$ As $|\alpha_{n+k}|>1/6$ and
$N(n)$ is increasing, it follows that
$$|f_{n+k+1}||\alpha_{n+k}-\alpha_{n+k+1}|\le
M(13|\alpha_{n+k}|+|\alpha_{n+k+1}|)a^{N(n+k)}.$$ Hence
$$|\alpha_{n+k+1}|\left(1-\frac{Ma^{N(n+k)}}{|f_{n+k+1}|}\right)\le
|\alpha_{n+k}|\left(1+13\,\frac{Ma^{N(n+k)}}{|f_{n+k+1}|}\right).$$
So, by (\ref{E:ass}) and (\ref{E:ind}),
$|\alpha_{n+k+1}|<4|\alpha_{n+k}|$.
\par Thus
$$|\alpha_{n+k}-\alpha_{n+k+1}|\le
17M2^{2n+6}a^{N(n+k)}2^{n+2}6^k\le2^{3n+3k+13}Ma^{N(n+k)},$$ and
by (\ref{E:ass})
$$|\alpha_{n+k}-\alpha_{n+k+1}|\le2^{-3k-3}<6^{-k-1}-6^{-k-2}.$$
Using the bounds for $|\alpha_{n+k}|$, this yields the desired
estimates for $|\alpha_{n+k+1}|$.
\par The inductive proof of the inequalities (\ref{E:ind}) is now
concluded. Moreover, we have shown that
$$|\alpha_m-\alpha_{m+1}|\le2^{3m+13}Ma^{N(m)},$$
for all $m\ge n$. This implies that
$\alpha_m\rightarrow\alpha\in{\Bbb C}$, and for $m\geq n$
$$|\alpha-\alpha_m|\le2^{13}M\sum_{j=m}^\infty2^{3j}a^{N(j)}\le
2^{13}Ma^{N(m)/2}\sum_{j=m}^\infty2^{3j}a^{N(j)/2}.$$ Hence
\begin{equation}\label{E:est}
|\alpha-\alpha_m|\le Ba^{N(m)/2},\;
B=2^{13}M\sum_{j=0}^\infty2^{3j}a^{N(j)/2}.
\end{equation}
\par Let $Q(z)=\alpha z-1$. Lemma \ref{L:limsup} implies that
$|\alpha|\geq1$. If $|\alpha|>1$ then by Lemma \ref{L:liminf} $f$
is a polynomial, which is in contradiction to $\rho=1$. Thus
$|\alpha|=1$. We let $$P(z)=Q(z)f(z)=\sum_{k\ge0}c_kz^k.$$ Note
that
$$P(z)-P_m(z)=Q_m(z)f(z)-P_m(z)+(\alpha-\alpha_m)zf(z).$$
It follows, using (\ref{E:mainest}), (\ref{E:ass}), (\ref{E:ind})
and (\ref{E:est}), that
\begin{eqnarray*}|c_{m+1}|&\le&|\alpha_mf_m-f_{m+1}|+
|\alpha-\alpha_m||f_m|\\&\leq&M\left(2^{2n+6}+1\right)a^{N(m)}+
2^mBa^{N(m)/2},\end{eqnarray*} for all $m\ge n$. This implies that
$|c_m|^{1/m}\rightarrow0$, hence $P$ is an entire function.
\par Observe that $Q_m(z)P(z)-(\alpha z-1)P_m(z)$ has $N(m)$
zeros in $\overline\Delta_r$. Since $P$ is entire, it follows by
Lemma \ref{L:liminf} that $P$ is in fact a polynomial. So $f=P/Q$,
and $Q$ does not divide $P$ since $f$ is not entire. This finishes
the proof. $\Box$

\par Theorem \ref{T:rational} has the following corollary, which
is proved exactly as Corollary \ref{C:Lagr}.
\begin{Corollary}\label{C:rational} Let $\{n_k\}_{k\geq0}$ be an
increasing sequence of natural numbers such that $n_{k+1}/n_k\leq
C$ for some constant $C$.  Let $f\in O(\Delta)$ and
$K=\overline\Delta_r$, where $r<1$. If
$\lim_{k\to\infty}N_K(n_k,1)/n_k=\infty$, then either $f$ is
entire or $f=P/Q$, where $P,\,Q$ are polynomials, $\deg Q=1$ and
$Q$ does not divide $P$.\end{Corollary}

\par We conclude this section with a remark about Pad\'e
overinterpolation. Let $f$ be a germ of a holomorphic function at
the origin. A rational function $R\in{\mathcal R}_{nm}$ is called
a Pad\'e interpolator (or Pad\'e approximant) of type $(m,n)$ of
$f$ if $f-R$ has a zero of the highest possible order at the
origin, i.e. of order $N_K(n,m)$, where $K=\{0\}$. We prove the
following simple fact about overinterpolation in the $m$-th row of
the Pad\'e table.
\begin{Proposition}\label{T:Pade} Let $f$ be a holomorphic germ at
the origin and $m\in{\Bbb N}$. If, for all $n\ge k$, there exist
functions $R_n\in{\mathcal R}_{nm}$ so that $f-R_n$ vanishes to
order at least $n+m+2$ at the origin, then $f\in{\mathcal
R}_{km}$.
\end{Proposition}
\begin{pf} Let us write $R_n=P_n/Q_n$, where $P_n\in{\mathcal P}_n$
and $Q_n\in{\mathcal P}_m$, $Q_n\neq0$. For $n\geq k$ the function
$R_n-R_{n+1}$ vanishes to order at least $n+m+2$ at the origin.
Since $\deg(P_nQ_{n+1}-P_{n+1}Q_n)\leq n+m+1$, this implies
$R_n=R_{n+1}=R\in{\mathcal R}_{km}$, for $n\geq k$. It follows
that $f=R$.\end{pf}

\section{Overinterpolation and overconvergence}\label{S:oao}
\par Throughout this section we assume that $f\in O(\overline\Delta)$
and that $0<r<1$ is fixed. For a compact set $E\subset{\Bbb C}$
and a continuous complex-valued function $g$ on $E$, we denote by
$\|g\|_E$ the uniform norm of $g$ on $E$.
\par The following theorem shows that, in the presence of
overinterpolation, the functions $R_{nm}$ quickly approximate $f$
on some circle $S_t=\{z\in{\Bbb C}:\,|z|=t\}$.
\begin{Theorem}\label{T:cl} Let $m(n)\in{\Bbb N}$, and
$d_n>0$ be so that $\sum d_n$ converges. Suppose that for all $n$
there are polynomials $P_n\in{\mathcal P}_n$ and
$Q_{m(n)}\in{\mathcal P}_{m(n)}$, $Q_{m(n)}\ne0$, so that the
function $Q_{m(n)}f-P_n$ has $N(n)$ zeros in $\overline\Delta_r$.
There exist positive constants $b<1$, $c$, depending only on $r$,
and $t\in[r,(1+r)/2]$, such that
$$\|f-R_n\|_{S_t}\le M\left(\frac{c}{d_n}\right)^{m(n)}b^{N(n)}$$
holds for all $n$ sufficiently large, where $R_n=P_n/Q_{m(n)}$ and
$M=M(1,f)$.
\end{Theorem}
\begin{pf} We may assume that $M(r/2,Q_{m(n)})=1$. Following
\cite{CP2}, we define the $n$-th diameter of a set $G\subset{\Bbb
C}$ by
$$\diam_n(G)=\inf\left\{r_1+\dots+r_k:\,k\leq
n,\;G\subset\bigcup_{j=1}^kC_j(r_j)\right\},$$ where $C_j(r_j)$
are closed disks of radii $r_j>0$. If $H_n(z)=Q_{m(n)}(rz/2)$ then
by Lemma 3.3 from \cite{CP2}, for every $0<h\leq 1/(8e)$, the
$n$-th diameter of the set
$$G'=\left\{z\in{\Bbb C}:\,|H_n(z)|\leq
\left(\frac{hr^2|z|}{(1+r)^2}\right)^{m(n)}, 2\leq|z|\leq
\frac{1+r}r\right\}$$ does not exceed $36eh$. Hence the $n$-th
diameter of the set
$$G=\left\{z\in{\Bbb C}:\,|Q_{m(n)}(z)|\leq
\left(\frac{2hr|z|}{(1+r)^2}\right)^{m(n)}, r\leq|z|\leq
\frac{1+r}2\right\}$$ does not exceed $18ehr$. This means that the
measure of the set
$$F_n=\left\{t\in\left[r,\frac{1+r}{2}\right]:\,
|Q_{m(n)}(z)|\ge\left(\frac{2hr|z|}{(1+r)^2}\right)^{m(n)},
\;\forall\,z\in S_t\right\}$$ is at least $(1-r)/2-36ehr$.
\par Since $M(r/2,Q_{m(n)})=1$, the classical Bernstein-Walsh
inequality implies that $$M(1,Q_{m(n)})\leq(2/r)^{m(n)}.$$ If
$t\in F_n$ then by (\ref{e:mfe0}) we have
$$M(t,Q_{m(n)}f-P_n)\leq M\left(\frac{2}{r}\right)^{m(n)}
\left(\frac3{1-t}\right)^{m(n)+2}a^{N(n)}_1,$$
where
$$a_1=a_1(t)=\frac{12t^2+6t}{13t^2+4t+1}<1.$$ The function $a_1(t)$ is
increasing on $[0,1]$ and, therefore, it does not exceed
$$b=a_1((1+r)/2)$$ on $F_n$. Hence for $t\in F_n$ we have
$$\|f-R_n\|_{S_t}\le
\frac{9M}{(1-t)^2}\left(\frac{3(1+r)^2}{hr^2t(1-t)}\right)^{m(n)}
b^{N(n)}\le M\left(\frac{c_1}h\right)^{m(n)}b^{N(n)},$$ where
$c_1$ is a constant depending only on $r$.
\par If we let $h=d_n/(36er)$ then the measure of
$F_n$ is at least $(1-r)/2-d_n$ and for $t\in F_n$ we have
$$\|f-R_n\|_{S_t}\le M\left(\frac{c}{d_n}\right)^{m(n)}b^{N(n)},$$
where $c=36erc_1$. Since $\sum d_n<\infty$ there is $n_0$ such
that the set $F=\bigcap_{n=n_0}^\infty F_n$ is not empty. If $t\in
F$ then the conclusion of the theorem holds for $t$ and for all
$n\ge n_0$.\end{pf}
\par If $g$ is a continuous function on
a compact set $E\subset{\Bbb C}$, we let
$$\rho(n,m)=\inf\|g-R\|_E,$$ where the infimum is taken over
all $R\in{\mathcal R}_{nm}$. We say that rational functions {\it
overconverge} to $g$ on $E$ if
$$\lim_{n\to\infty}\rho(n,m(n))^{1/n}=0,$$ for some sequence
$m(n)\in{\Bbb N}$.
\par The following corollary shows that, under suitable conditions,
overinterpolation implies overconvergence.
\begin{Corollary}\label{C:oioc} Under the assumptions of Theorem
\ref{T:cl}, suppose that there is a sequence $\{a_n\}$ of positive
numbers converging to 0 such that
$$\sum_{n=1}^\infty\frac{b^{N(n)/m(n)}}{a_n^{n/m(n)}}<\infty.$$
Then there exists $t\in[r,(1+r)/2]$ for which
$$\lim_{n\to\infty}\|f-R_n\|^{1/n}_{S_t}=0.$$
\end{Corollary}
\begin{pf} For the proof, take $d_n=c\,b^{N(n)/m(n)}/a_n^{n/m(n)}$.
\end{pf}
\par The fact that overinterpolation implies overconvergence
allows us to use results of Gonchar and Chirka to prove other
results about overinterpolation. Let us first recall some
definitions from \cite{Ch}. The class ${\mathcal R}_{n,(m)}$
consists of all rational functions of degree at most $n$ and with
at most $m$ geometrically distinct poles. The class ${\mathcal
A}^0_m$ consists of all functions meromorphic on ${\Bbb P}^1$
except for at most $m$ singularities of finite order. This means
that for every singular point $a$ there is a number $p$ such that
$|f(z)|<\exp(1/|z-a|^p)$ near $a$.
\begin{Theorem}\label{T:gc} Let $m\ge0$ be an integer. If
for all $n$ there are functions $R_n\in{\mathcal R}_{nm}$ such
that $f-R_n$ has $N(n)$ zeros in $\overline\Delta_r$, where
$N(n)/n\rightarrow\infty$, then $f$ extends to a meromorphic
function on ${\Bbb C}$ with at most $m$ poles.
\par If the functions $R_n\in{\mathcal R}_{n,(m)}$ and
$$\liminf_{n\to\infty}\frac{N(n)}{n\log n}>-\frac1{\log b}\;,$$
where $b$ is the constant from Theorem \ref{T:cl}, then $f$ has an
extension in ${\mathcal A}^0_m$.
\end{Theorem}
\begin{pf} To prove the first statement, we take a number $\alpha$ such
that $b<\alpha<1$ and let $a_n=\alpha^{N(n)/n}$. By Corollary
\ref{C:oioc}, there is $t\in[r,(1+r)/2]$ for which
$$\lim_{n\to\infty}\|f-R_n\|^{1/n}_{S_t}=0.$$ By Theorem 1 from
\cite{Go}, $f$ extends to a meromorphic function to ${\Bbb C}$
with at most $m$ poles.
\par A result of Chirka and Gonchar (see \cite[Theorem 1]{Ch})
states that if $f$ is analytic in a neighborhood of a compact set
$E\subset{\Bbb C}$ of positive capacity then $f$ has an extension
in ${\mathcal A}^0_m$ if and only if there are a sequence of
rational functions $R_n\in {\mathcal R}_{n,(m)}$ and a number
$\lambda>0$ such that
$$\|f-R_n\|^{1/n}_E<\frac1{n^\lambda}$$ for all $n$ sufficiently
large. (The theorem is stated for $E=\overline\Delta_s$, but see
the note after the statement.)
\par Take numbers $\alpha$ and $\lambda$ such that
$$\liminf_{n\to\infty}\frac{N(n)}{n\log n}>\alpha>-\frac1{\log b}
\;,\;\;0<\lambda<-1-\alpha\log b.$$ Let $$d_n=cb^{\alpha\log
n}n^\lambda=cn^{\alpha\log b+\lambda},$$ where $c$ is the constant
from Theorem \ref{T:cl}. Then $\sum d_n<\infty$ and by Theorem
\ref{T:cl} there is $t\in[r,(1+r)/2]$ such that
$$\|f-R_n\|^{1/n}_{S_t}\le M^{1/n}\frac{c}{d_n}\,b^{N(n)/n}=
M^{1/n}\frac{b^{N(n)/n-\alpha\log n}}{n^\lambda}<
\frac1{n^\lambda}\;,$$ for all $n$ sufficiently large. Now the
second statement of the theorem follows from the result of Chirka
and Gonchar mentioned above.
\end{pf}

\section{Interpolation by algebraic functions}\label{S:obaf}
\begin{Proposition}\label{P:revBezout} Let $S$ be an infinite set
in ${\Bbb C}^2$ with the following property: There exist positive
constants $A\ge1$ and $\alpha<2$ such that $$|S\cap X|\leq A(\deg
X)^\alpha,$$ for any algebraic curve $X\subset{\Bbb C}^2$ not
containing $S$. Then $\alpha\geq1$ and $S$ is contained in an
irreducible algebraic curve of degree at most
$(2A)^{1/(2-\alpha)}$. Moreover,
$$|S\cap X|\leq(2A)^{1/(2-\alpha)}\,\deg X,$$ for any
algebraic curve $X\subset{\Bbb C}^2$ not containing $S$.
\end{Proposition}
\begin{pf} Suppose $\alpha<1$. Assume that
$\{z_1,\dots,z_n\}\subseteq S$, where $n\geq2$, and let $L_j$,
$1\leq j<n$, be a complex line passing through $z_j$ and not
containing $z_n$. If $X=L_1\cup\dots\cup L_{n-1}$, then $X$ does
not contain $S$, hence
$$n-1\leq|S\cap X|\leq A(n-1)^\alpha.$$ Thus
$|S|\leq1+A^{1/(1-\alpha)}$, which is a contradiction.
\par Let $k$ denote the greatest integer in $x=(2A)^{1/(2-\alpha)}$.
Then $$2Ak^\alpha= x^{2-\alpha}k^\alpha<(k+1)^2\leq k^2+3k.$$ Note
that the dimension of the space of polynomials in ${\Bbb C}^2$ of
degree at most $n$ is $(n+1)(n+2)/2$. Therefore there exists a
curve $C$ of degree at most $k$ so that
\begin{equation}\label{E:est2}|S\cap C|\geq(k^2+3k)/2>Ak^\alpha.
\end{equation} It follows that $S\subseteq C$. Assume that
$C=C_1\cup\dots\cup C_m$, where $C_j$ is an irreducible algebraic
curve of degree $k_j$, $k_1+\dots+k_m\leq k$. If no curve $C_j$
contains $S$ then, since $\alpha\geq1$, $$|S\cap
C|\leq\sum_{j=1}^m|S\cap C_j|\leq A\sum_{j=1}^mk_j^\alpha\leq
Ak^\alpha,$$ which contradicts (\ref{E:est2}). We conclude that
$S$ is contained in an irreducible curve $\Gamma$ of degree at
most $k$. Hence by Bezout's theorem,
$$|S\cap X|\leq|\Gamma\cap X|\leq(2A)^{1/(2-\alpha)}\,\deg X,$$
for any algebraic curve $X\subset{\Bbb C}^2$ not containing $S$.
\end{pf}


\begin{thebibliography}{XXXX}
\bibitem[Ch]{Ch} E. M. Chirka, {\em Rational approximations of holomorphic
functions with singularities of finite order}, Mat. Sb. (N.S.)
{\bf 100(142)} (1976), 137--155.
\bibitem[CP1]{CP1} D. Coman and E. A. Poletsky, {\em Bernstein-Walsh
inequalities and the exponential curve in ${\Bbb C}^2$}, Proc.
Amer. Math. Soc. {\bf 131} (2003), 879--887.
\bibitem[CP2]{CP2} D. Coman and E. A. Poletsky, {\em Measures of
transcendency for entire functions}, Mich. Math. J. {\bf 51}
(2003), 575-591.
\bibitem[CP3]{CP3} D. Coman and E. A. Poletsky, {\em Transcendence
measures and algebraic growth of entire functions}, preprint 2004.
\bibitem[G]{G} D. Gaier, {\em Lectures on Complex Approximation}, Birkh\"auser,
1987.
\bibitem[Go]{Go} A. A. Gonchar, {\em On a theorem of Saff},
Mat. Sb. (N.S.) {\bf 94(136)} (1974), 152--157.
\end{thebibliography}
\end{document}